\documentclass[11pt]{article}

\usepackage[usenames]{color}
\usepackage{amssymb,amsmath,latexsym}
\usepackage{amsthm}
\usepackage[mathscr]{eucal}
\usepackage{setspace}
\usepackage{geometry}
 \usepackage{graphicx}
 \usepackage{float}
\usepackage{url}
    \usepackage{soul}
\usepackage{cancel}

\newtheorem{lemma}{Lemma}
\newtheorem{theorem}{Theorem}
\newtheorem{corollary}{Corollary}

\newtheorem{example}{Example}
\newtheorem{proposition}{Proposition}

\newtheorem{remark}{Remark}

\newcommand{\noi}{\noindent}
\newcommand{\beqs}{\begin{equation*}}
\newcommand{\eeqs}{\end{equation*}}
\newcommand{\beq}{\begin{equation}}
\newcommand{\eeq}{\end{equation}}
\newcommand{\beqys}{\begin{eqnarray*}}
\newcommand{\eeqys}{\end{eqnarray*}}
\newcommand{\beqy}{\begin{eqnarray}}
\newcommand{\eeqy}{\end{eqnarray}}

\newcommand{\N}{I\!\!N}

\newcommand{\PP}{I\!\!P}

\newcommand{\E}{I\!\!E}
\newcommand{\R}{I\! \! R}

\newcommand{\pos}{\in(0,\infty)}

\newcommand{\ra}{\rightarrow}
\newcommand{\Ra}{\Rightarrow}

\newcommand{\sig}[2]{\sum_{#1}^{#2}}
\newcommand{\inte}[2]{\int_{#1}^{#2}}

\newcommand{\wt}{\widetilde}

\newcommand{\var}{\operatorname{var}}

\mathchardef\mhyphen="2D
\newcommand{\bp}{\begin{pmatrix}}
\newcommand{\ep}{\end{pmatrix}}
\newcommand{\black}[1]{{\textcolor{black}{#1}}{}}

\begin{document}

\title{Stationary Distribution Convergence of the Offered Waiting Processes 
in Heavy Traffic under General Patience Time Scaling}

\author{Chihoon Lee\footnote{E-mail: Chihoon.Lee@stevens.edu}\\ School of Business\\ Stevens Institute of Technology\\ Hoboken, NJ 07030 \\
\and Amy R. Ward\footnote{E-mail: amy.ward@chicagobooth.edu}\\Booth School of Business\\The University of Chicago\\ Chicago, IL 60637 \\
\and Heng-Qing Ye\footnote{E-mail: lgtyehq@polyu.edu.hk}\\Dept. of Logistics and Maritime Studies\\Hong Kong Polytechnic University\\ Hong Kong}
\date{\today}\maketitle

\begin{abstract}
\noi  
We study a sequence of single server queues with 
customer abandonment ($GI/GI/1+GI$) under heavy traffic. 
The patience time distributions vary with the sequence, which allows for a wider scope of applications. 
It is known (\cite{ReedWard2008,LeeWeerasinghe2011}) that the sequence of scaled offered waiting time processes converges weakly to a reflecting diffusion process with non-linear drift, as the traffic intensity approaches one.
In this paper, we further show that the sequence of stationary distributions and moments of the offered waiting times, with diffusion scaling, converge to those of the limit diffusion process. This justifies the stationary performance of the diffusion limit as a valid approximation for the stationary performance of the $GI/GI/1+GI$ queue. Consequently, we also derive the approximation for the abandonment probability for the $GI/GI/1+GI$ queue in the stationary state.
\end{abstract}

\noi {\it Keywords: Customer Abandonment; Heavy Traffic; Stationary Distribution Convergence}

\section{Introduction}

In a recent paper,  Lee {\it et al.}~\cite{LeeWardYe-QueSyst}, we have studied a sequence of single server systems with abandonment ($GI/GI/1+GI$). In each system, the customers arrive following a renewal process and have a general service time distribution. 
We have shown that under appropriate conditions the stationary distributions of suitably scaled offered waiting time processes converge to the stationary distribution of the associated reflected Ornstein-Uhlenbeck process  in the positive real line, as the traffic intensity approaches one (that is, as the heavy traffic scaling parameter $n$ approaches infinity).
One of the key assumptions made in the analysis is that the patience times have remained the same and \emph{unscaled} for all the systems, and consequently it leads to a diffusion approximation that depends only on the behavior of the patience time distribution $F(\cdot)$ at 0. 
Such diffusion approximations, with a linear drift $F'(0)x$ accounting for the customer abandonment, have at least two drawbacks: (Applicability) There is a wide family of distributions where $F'(0) = 0$ or $\infty$, e.g., Gamma distributions with the shape parameter $\alpha \neq 1$; (Performance) If the underlying patience time distribution is of increasing (or decreasing) hazard rate (i.e., the longer one has waited for service, the more likely to abandon), then the diffusion approximation with a linear drift $F'(0)x$ invariably over-estimates (or under-estimates) the steady-state offered waiting time.

Our objective in this paper is to establish the stationary distribution convergence of offered waiting time in the heavy traffic limit whilst the underlying patience times are allowed to vary and are suitably scaled, ensuring the entire distributional information is contained in the heavy traffic diffusion limit. 
This relaxes the assumption of uniform abandonment time distribution in our previous work and greatly enhances the applicability of our results. 
As a result, we derive the approximation for the abandonment probability for the $GI/GI/1+GI$ queue in the stationary state.
It is also interesting to note that we are able to bypass the hydrodynamic scaling approach in Lee {\it et al.}~\cite{LeeWardYe-QueSyst} and simplify the arguments based on more direct and simpler pathwise stability estimates.     



Our result draws on past work that has developed heavy traffic approximations for the $GI/GI/1+GI$ queue using the offered waiting time process.  The offered waiting time process, introduced in \cite{Baccelli84}, tracks the amount of time an infinitely patient customer must wait for service. 
Its heavy traffic limit when the patience time distribution is left unscaled is a reflected Ornstein-Uhlenbeck process (see Ward and Glynn \cite{WardGlynn2005}), and its heavy traffic limit when the patience time distribution is scaled through its hazard rate is a reflected nonlinear diffusion (see Reed and Ward \cite{ReedWard2008}). The work of Lee and Weerasinghe \cite{LeeWeerasinghe2011} incorporates both of these scaling scenarios in the general framework that can be satisfied by many other classes of patience time distributions. One key ingredient in the proof of their main result is the martingale functional central limit theorem, which helps to accommodate the more general assumptions (e.g., the customer arrival process has a state dependent intensity).    

Generally speaking, those results (viz. the process-level convergence) are not sufficient to conclude that the stationary distribution of the offered waiting time process converges, which is the key to approximating the performance measures such as the stationary abandonment probability and mean queue-length.  Those limits were conjectured in~\cite{ReedWard2008}, and shown through simulation to provide good approximations.  However, the proof of those limits was left as an open question. In this paper, we analyze the case when the patience time distribution is scaled, following the general framework in  \cite{LeeWeerasinghe2011}. 



There is a growing literature on the problem of stationary distribution convergence of service systems in heavy traffic. Broadly speaking, there are \emph{indirect} and \emph{direct} methods in approaching the problem. 
As for the \emph{indirect} method, the papers \cite{Gamarnik:2006} and \cite{BudhLee07} establish the validity of the heavy traffic stationary approximation (viz. interchange of limits) for a generalized Jackson network, without customer abandonment, using a Lyapunov function method.  
The main difficulty in extending their methodologies to the current model is 
the lack of the global Lipschitz continuity property of an associated regulator (Skorokhod) mapping that helped
convert the given moment bound of primitives (the inter-arrival and
service times) to the bound of the waiting time or queue-length processes. The known regulator mapping under customer abandonment is only \emph{locally} Lipschitz (that is, the Lipschitz constant depends on time parameter). 

As for the more \emph{direct} methods, there are the MGF--BAR (moment generating function--basic adjoint relationship) method and the generator comparison method (also known as Stein's method). For the latter, the authors of \cite{HuangGurvich18} study the Poisson arrival case (i.e., $M/GI/1+GI$ queue) and show the associated Brownian model is accurate uniformly over a family of patience time distributions and universally in the heavy-traffic regime. Owing to the Poisson arrivals, it is enough to consider a one-dimensional process with a simple generator, whereas with general arrival processes, one needs to consider a two-dimensional process (tracking, e.g., the residual arrival times) and correspondingly more  complicated generator. For the former approach, the authors of \cite{BravermanDaiMiyazawa2017heavy} consider a generalized Jackson network, without customer abandonment, and work directly with the BAR, an integral equation that characterizes the stationary distribution of a Markov process, and establish the convergence of the MGFs of the pre-limit stationary distributions.  When extending their methodology to the current model, a unique challenge arises in establishing a convergence rate (in terms of the scaling parameter) of the tail estimates for the stationary abandonment probability. Lastly, there is the so-called Drift method; see \cite{eryilmaz2012asymptotically} and also \cite{hurtado2020transform} for its connection to the MGF method. The Drift method uses polynomial test functions in an inductive manner by setting to zero the drift of the test function (i.e., equating the expected value of the test function in two different time steps).   


In view of the above-mentioned technical challenges, we adopt the approach in the studies of \cite{YeYao12-interchange, YeYao18}, which extend the works of \cite{Gamarnik:2006} and \cite{BudhLee07} (the interchange of limits) to a wider range
of stochastic processing networks, e.g., the multiclass queueing network
and the resource-sharing network.
In their studies, they relax the requirement of the aforementioned Lipschitz continuity by establishing the so-called uniform stability and the uniform moment bound. Further extending their approach to the {\it nonlinear} dynamic complementarity problem, which is relevant to the $GI/GI/1+GI$ model,  Lee {\it et al.}~\cite{LeeWardYe-QueSyst} establishes the stationary distribution and moment convergences when the patience times are left unscaled.  
The key proof of Lee {\it et al.}~\cite{LeeWardYe-QueSyst} is to establish a kind of uniform stability property (cf. Lemmas \ref{lem-10a}-\ref{lem-10cd}) by applying Bramson's hydrodynamic scaling approach and its variation. (By the hydrodynamic approach, the $n$th diffusion-scaled process breaks into many pieces of fluid-scaled processes, with each piece covering a period of $1/\sqrt n$ in the diffusion-scaled process.)

The contributions of the current paper are: (a) identifying the stability conditions of the pre-limit system (Theorem \ref{stability}) under general scaling assumption; (b) establishing the convergence of stationary distributions and moments of offered waiting times (Theorem~\ref{convstat}) and offering simplified proofs based on pathwise stability estimates, replacing the more involved hydrodynamic scaling approach taken in \cite{LeeWardYe-QueSyst}, {and (c) establishing the convergence of the scaled stationary abandonment probabilities  (Corollary~\ref{aban-prob}).} 

The remainder of this paper is organized as follows. 
In Section \ref{subsection:modelassumption}, we set up the model assumptions and recall the known process-level convergence results for the $GI/GI/1+GI$ queue under general patience time distribution scaling.  In Section \ref{section:SSmoments}, we state our mains results on the convergence of the stationary distribution of the offered waiting time process and its moments.    
In Section \ref{section:UniformMoment}, we present key moment bounds of the scaled state processes that  are  uniform in the heavy traffic scaling parameter ($n$). Lastly, in Section \ref{section:LemmaProof}, we provide the proofs of key lemmas.  


{\bf Notation and Terminology.} Use the symbol ``$\equiv$'' to stand for equality by definition. The set of positive integers is denoted by $\N$ and denote $\N_0\equiv \N\cup\{0\}$.  Let $\R$ represent the real numbers $(-\infty, \infty)$ and $\R_+$ the non-negative real line $[0,\infty)$.  For $x,y \in \R$, $x \vee y \equiv \max\{x,y\}$ and $x \wedge y \equiv \min\{x,y\}$.   
Let $D(\R)\equiv D(\R_+, \R) $ be the space of right-continuous functions $f:\R_+\ra\R$ with left limits, endowed with the Skorokhod $J_1$-topology (see, for example,~\cite{Billingsley:1999}). Lastly, the symbol ``$\Ra$'' stands for the weak convergence; we make this explicit for stochastic processes in $D(\R)$, otherwise, it is used for weak convergence for a sequence of random variables. 


\section{The Model and \black{Known} Results} \label{subsection:modelassumption}
We consider a sequence of single server systems having FIFO service with abandonment indexed by $n\in \N$, and by  convention, we use superscript $n$ for any processes or quantities associated with the $n$-th system. For $n\in\N$, consider three independent i.i.d. sequences of nonnegative random variables $\{ u_i^n, i \geq 2\}$, $\{v_i^n,i\geq 1\}$, $\{d_i^n,i\geq 1\}$, that are representing inter-arrival times, service times, and patience times, respectively, and are defined on a common probability space $(\Omega, \mathcal{F}, \PP)$. 

We first consider the sequences $\{ u_i^n, i \geq 2\}$, $\{v_i^n,i\geq 1\}$ which are built from independent i.i.d. sequences of random variables with unit mean $\{ u_i, i \geq 2\}$, $\{v_i,i\geq 1\}$ in the following way. At time 0, the previous arrival to the system occurred at time $t_0^n <0$, so that $|t_0^n|$
represents the time elapsed since the last arrival in the $n$-th system.  We let $u_1$ be the random variable representing the remaining time conditioned on $|t_0^n|$ time units having passed; that is,
\[
\PP(u_1 >x) = \PP(u_2>x | u_2 > |t_0^n| ).
\]
Given positive sequences $\{\lambda^n\}$ and $\{\mu^n\}$, the $i$-th arrival to the $n$-th system occurs at time
\[
t_i^n \equiv \sum_{j=1}^i u_j^n\,, \quad  u_j^n\equiv\frac{u_j}{\lambda^n}\,, 
\]
{and} 
has service time
\[
    v_i^n \equiv \frac{v_i}{\mu^n}\,, 
\]
and abandons without receiving service if processing does not begin by time $t_i^n + {d_i^n}$. 

We assume the following conditions.
\begin{enumerate}
\item[($\mathbb A$1)] {For some $p\in (2,\infty)$, {$\E[u_2^{p}+v_2^{p}]<\infty$.}} 
\end{enumerate} 
The arrival and service rates in the $n$-th system, $\lambda^n$ and $\mu^n$, respectively, satisfy the following heavy traffic assumption:
\begin{enumerate}
\item[($\mathbb A$2)] $\lambda^n\equiv n\lambda$, $\lim_{n \rightarrow \infty} \frac{\mu^n}{n} = \lambda \in (0,\infty)$ and $\lim_{n \rightarrow \infty} \sqrt{n} \left(\lambda - \frac{\mu^n}{n} \right) = \theta \in \R.$
\end{enumerate}

Next, we consider the following assumption on the patience time distributions $\{d_i^n,i\geq 1\}$ as studied in Lee and Weerasinghe \cite{LeeWeerasinghe2011}. As a consequence, the drift coefficient of the limiting diffusion is influenced by the sequence of patience time distributions in a non-linear fashion. 
\begin{enumerate}
\item[$(\mathbb A3)$] 
Let $F^n(\cdot)$ be the right continuous  patience time distribution function of the i.i.d. sequence $(d_i^n)_{i\geq1}$.  Assume that $F^n(0)=0$ and there is a non-negative continuously differentiable  function $H(\cdot)$  such that for each $K>0$, $$\lim\limits_{n\ra\infty}\sup\limits_{x\in[0,K]}\left|\sqrt n F^n\left(\frac{x}{\sqrt n}\right)-H(x)\right|=0.$$\end{enumerate} 

\begin{example}\label{eg1}  We provide some examples of patience time distribution functions  $\{F^n\}$ satisfying $(\mathbb A3)$; cf.  Lee and Weerasinghe \cite{LeeWeerasinghe2011} and Huang {\it et al.}~\cite{huang2016unified}.  
\begin{enumerate} 
\item Take $F^n = F$ for all $n$, where $F(\cdot)$ is some distribution function, differentiable with a bounded derivative on $[0,\delta]$ for some $\delta>0$.  Hence, $H(x)=F'(0)x$ in this case.  This corresponds to the result in Ward and Glynn \cite{WardGlynn2005}. 
\item Take $F^n(x)=1-\exp(-\inte{0}{x}h(\sqrt n u)du)$ for $x\geq0$, where $h(\cdot)$ is a continuous hazard rate function.  In this case, $H(x)=\inte{0}{x}h(u)du$ and it satisfies $(\mathbb A3)$.  Indeed, for any general sequence $\{F^n\}$, if $F^{n}(\frac{x}{\sqrt n})$ converges to a non-negative function $h(x)$ uniformly on compact sets, then $\{F^n\}$ satisfies $(\mathbb A3)$ with the limiting function $H(x)=\inte{0}{x}h(u)du$. This corresponds to the result in Reed and Ward \cite{ReedWard2008}.  \item 
Take any non-negative, non-decreasing, continuously differentiable function $H(\cdot)$  which satisfies $H(0)=0$ and $H(\infty)=\infty$.  Let $F^n(x)=\frac{1}{\sqrt n} \min\{H(\sqrt n x),\sqrt n\}$ for all $x\geq0$.  Then, for each $n\geq1$, $F^n$ is a continuous probability distribution function and the sequence of distribution functions $\{F^n\}$ satisfies$(\mathbb A3)$ with limiting function $H(\cdot)$.  
\item {More examples such as mixture of hazard-rate scaling and no scaling, and delayed hazard-rate scaling are presented in Huang {\it et al.}~\cite{huang2016unified}.}  \end{enumerate}
\end{example}

\noi \textbf{The Offered Waiting Time Process}


The \emph{offered waiting time} process $ \{V^n(t):t\geq0\}$ tracks the amount of time an incoming customer has to wait for service: \beq \label{2.1.1}  V^n(t)=  \black{V^n(0) } +  \sig{j=1}{A^n(t)} v^n_j\mathbf{1}_{[V^n(t^n_j-)<d_j^n]} - \inte{0}{t}\mathbf{1}_{[V^n(s)>0]}ds\geq0. \eeq
Here, \beq\label{254} A^n(t) \equiv \max\{ i \in \N_0 : t_i^n \leq t \} \eeq
 is a delayed renewal process whenever $|t_0^n| >0$.

{\noi {\bf Reflected Diffusion Approximation with General Patience Time Distribution}

Consider the one-dimensional reflected diffusion process $V\equiv \{V(t): t\geq 0\}$:
 \beq\label{3.58}
 \begin{array}{l}
V(t)=V(0) +  {\sigma} W(t)+ {\frac{\theta}{\lambda}} t- {\inte{0}{t}H(V(s))ds} + L(t) \geq 0 \\
 \mbox{subject to:  } L \mbox{ is non-decreasing, has } L(0)=0 \mbox{ and } \int_{0}^{\infty}V(s)dL(s)=0,
 \end{array} \eeq
where $\{W(t):t\geq0\}$ denotes a one-dimensional standard Brownian motion, and the infinitesimal variance parameter is
 $\sigma^2  \equiv \lambda^{-1}( \var(u_2) + \var(v_2)).$


{The following weak convergence result is a simple modification of Theorem 4.10 of Lee and Weerasinghe \cite{LeeWeerasinghe2011} (see also Reed and Ward \cite{ReedWard2008} and  Lee {\it et al.}~\cite{LeeWardYe-QueSyst}). 
\begin{proposition}\label{pc}
Assuming $\sqrt{n} V^n(0) \Rightarrow V(0)$ as $n\ra\infty$, we have
  \begin{equation} \label{eq:Vapprox}
 \sqrt{n}V^n \Ra V \,\,\mbox{ in }\,\, D(\R) \,\,\mbox{ as } \,\,  n \ra \infty.
 \end{equation}
 \end{proposition}
The weak convergence (\ref{eq:Vapprox}) motivates approximating the scaled stationary distributions for $V^n$, and its moments, using the stationary distribution of $V$, and its moments. First of all, we impose the following well-known stability condition for a one-dimensional diffusion $V$ (cf. \cite{Echeverria} and Proposition 6.1(i) in~\cite{ReedWard2008}):
\begin{enumerate}
\item[($\mathbb A4$)]  $\displaystyle\lim_{x\rightarrow \infty} H(x) > \frac{\theta}{\lambda}$\,.
\end{enumerate}
We note that the above stability condition is trivially satisfied for the two examples described earlier, right below $(\mathbb A3)$. That is, both $H(x)=F'(0)x$ and $H(x)=\inte{0}{x} h(u)du$ converge to $\infty$ as $x\ra\infty$.
Identical arguments as in the proof of Proposition 6.1(i) in~\cite{ReedWard2008} show that, under  ($\mathbb A4$),
\beq\label{lim-int}
V(t) \Rightarrow V(\infty) \,\mbox{ as }\, t \rightarrow \infty
\eeq
for $V(\infty)$ a random variable having a density function \beq\label{eq:VstationaryDensity}f(x)=M\exp\left(\frac{2}{\sigma^2}\left(\frac{\theta}{\lambda} x-\inte{0}{x}H(s)ds\right)\right), \quad x\geq0,\eeq 
where $M\pos$ is such that $\inte{0}{\infty}f(x)dx=1$. 
In summary, the convergences \eqref{eq:Vapprox} and \eqref{lim-int} correspond to the limit when $n\rightarrow \infty$ is taken first and $t \rightarrow \infty$ is taken second, depicted graphically in Figure~\ref{120}.}  

\section{{The Stationary Distribution Existence and Convergence Results}} \label{section:SSmoments}



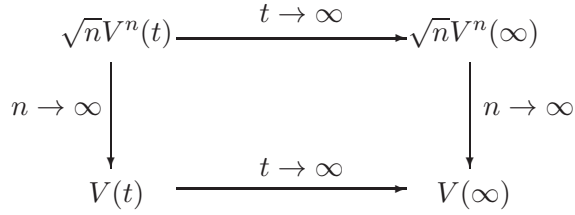
\begin{figure}[!t]
\centering
\unitlength 1mm 
\linethickness{0.4pt}
\ifx\plotpoint\undefined\newsavebox{\plotpoint}\fi 
\begin{picture}(75.75,31)(50,115)
\put(64,139.5){\makebox(0,0)[cc]{$\sqrt{n}V^n(t)$}}
\put(111.5,139.5){\makebox(0,0)[cc]{$\sqrt n V^n(\infty)$}}
\put(64,118){\makebox(0,0)[cc]{$V(t)$}}
\put(111.5,118){\makebox(0,0)[cc]{$V(\infty)$}}
\put(72,138.75){\vector(1,0){30.5}}
\put(72,118.75){\vector(1,0){30.5}}
\put(63.25,135.5){\vector(0,-1){14}}
\put(111,135.5){\vector(0,-1){14}} \put(83,141){$t\rightarrow
\infty$} \put(83,120.5){$t\rightarrow \infty$}
\put(50,128.5){${n\rightarrow\infty}$}
\put(112.75,128.5){${n\rightarrow\infty}$}
\end{picture}
\caption{\em A graphical representation of the limit interchange. } \label{120}
\end{figure}


Because of the {remaining} arrival time  (i.e., the forward recurrence time of the arrival process), the offered waiting time process $\{V^n(t):t\geq0\}$ alone is not Markovian. 
Defining the remaining arrival time
\[
\tau^n (t)\equiv t_{j+1}^n - t \,\mbox{ for }\, t \in [t_j^{{n}}, t_{j+1}^{{n}}), \,\,j \in \N_0,
\]
where {$\tau^n(0) =u_1$},  the process $\mathbb X^n \equiv \{(\tau^n(t), V^n(t)):t\geq0\}$ having state space $\mathbb S \equiv  \R_+ \times \R_+$ is strong Markov (cf. Problem 3.2 of Chapter X in \cite{Asmussen2003}). 
For $x=(\tau, v) \in \mathbb S$, define its norm $|x|$ as $|x|\equiv\tau+v$ and likewise define $|\mathbb X^n(t)| \equiv \tau^n(t) + V^n(t), \,\, t\geq0\,. $
We will assume the interarrival times are unbounded. Such an assumption has been frequently used in the literature 
to verify a \emph{petite} set requirement that implies positive Harris recurrence of a Markov process, cf. Proposition 4.8 in \cite{BramsonBook2008}. 

\begin{enumerate}
\item[($\mathbb A5$)] The i.i.d. interarrival times $\{u_i, i\geq2\}$ are unbounded, that is, $\PP(u_2\geq u)>0$ for any $u>0$. 
\end{enumerate}


\begin{theorem}\label{stability} \textsc{(Stationary Distribution Existence)} {Assume $(\mathbb A1)$--$(\mathbb A5)$.}
For any sufficiently large $n$,
there exists a unique stationary probability distribution for the Markov process ${\mathbb X}^n$.
\end{theorem}

Next, 
we consider the diffusion-scaled process:
\[  \widetilde{\mathbb X}^n \equiv (\wt\tau^n, \wt{V}^n), \mbox{ where } \wt{\tau}^n(t)\equiv \sqrt n \tau^n(t) ,\,\,\, \wt{V}^n(t)\equiv \sqrt n V^n(t)\,. \]
Notice 
the time is not scaled in the process $\wt{V}^n$ because the arrival and service rate parameters are scaled instead (from $(\mathbb A2)$, both $\lambda^n$ and $\mu^n$ are order $n$ quantities); therefore, scaling the state by $\sqrt{n}$ produces the traditional diffusion scaling. 

\begin{remark}\label{rem11}
{We note that Theorem~\ref{stability} requires sufficiently large $n$, while the result in Lee {\it et al.}~\cite{LeeWardYe-QueSyst} does not. This is because the proof of  Theorem~\ref{stability} is based on a weaker version of an estimate \eqref{2.1} in Proposition~\ref{moment1} below, which is $\lim_{|x|\ra\infty}\E\left[|\wt{\mathbb X}^n_x(t|x|)|^{{q}} \right]/{|x|^{{q}}}=0$ (there is no `supremum over $n$'), and we could establish this only for sufficiently large $n$. More precisely, we could prove Lemma \ref{lem-10a} below only for large enough (but fixed) $n$ (under which the negative drift of netput rate is warranted; see Remark \ref{rem22}).}
\end{remark}

\begin{theorem}\label{convstat} \textsc{(Stationary Convergence)} Assume $(\mathbb A1)$--$(\mathbb A5)$. For large enough $n$, let $\pi^n$ denote the stationary distribution of  $\widetilde{\mathbb X}^n$.

\begin{enumerate}
\item[(a)] \textsc{(Distribution)}    
Denote by $\pi^n_0$ the marginal distribution of $\pi^n$ on the {second} coordinate of {$\widetilde{\mathbb X}^n$,}  i.e., $\pi^n_0(A)= \pi^n(\R_+ \times A)$ for  $A\in\mathcal B(\R_+)$. Let $\wt{V}^n(\infty)$ be a random variable having distribution $\pi^n_0$ and also $V(\infty)$ a random variable having density \eqref{eq:VstationaryDensity}.
We have that $\widetilde{V}^n(\infty) \Rightarrow V(\infty)$ as $n \rightarrow \infty$.
\item[(b)] \textsc{(Moments)}
{For any  $m\in(0,p-1)$,} 
\begin{equation*} \E[(\wt V^n(\infty))^m] \ra \E[(V(\infty))^m] \,\mbox{ as }\, n \rightarrow \infty.  \end{equation*}
\end{enumerate}
\end{theorem}
 
As a consequence, we obtain the following result on the convergence of stationary abandonment probability. For $n\geq1$, define \[P^n_a\equiv \E[F^n(V^n(\infty))]= \E\left[ F^n\left(\frac{\wt V^n(\infty)}{\sqrt n}\right)\right],\] that is, a fraction of customers who abandon the system in stationarity. 
\begin{corollary}\label{aban-prob} \textsc{(Stationary Abandonment Probability Convergence)}
Assume $(\mathbb A1)$--$(\mathbb A5)$.  Assume further $\sqrt n F^n(x/\sqrt n)\leq C (1+x^m)$ for $m\in (0,p-1)$, where $C\pos$ is independent of $n$. Then, as $n\ra\infty$, \[\sqrt n P^n_a \ra \E[H(V(\infty))], \] where the function $H(\cdot)$  is as in $(\mathbb A3)$.
\end{corollary}
\begin{example}\label{eg2} We provide some examples of patience time distribution functions  $\{F^n\}$ satisfying the assumption  $\sqrt n F^n(x/\sqrt n)\leq C (1+x^m)$ in Corollary \ref{aban-prob} as well as  the scaling assumption $(\mathbb A3)$.
\begin{enumerate} 
\item Take $F^n =  F$ for all $n$ (as in Ward and Glynn \cite{WardGlynn2005}), where $F$ is differentiable with a derivative of polynomial growth satisfying $\sup_{y\in[0,x]} F'(y)\leq C(1+x^{m-1})$.
\item Take $F^n(x) = 1-\exp(-\inte{0}{x}h(\sqrt n u)du)$ (as in Reed and Ward \cite{ReedWard2008}), where $h(\cdot)$ is such that $\sup_{y\in[0,x]} h(y)\leq C(1+x^{m-1})$. 
\item For a general sequence $\{F^n\}$,  if $F^{n}(x/{\sqrt n})$ converges to a non-negative function $h(\cdot)$ uniformly on compact sets and $0\leq (F^n)'(x/\sqrt n)\leq C(1+x^{m-1})$ for some $C\pos$ independent of $n$, then $\{F^n\}$ satisfies the stated assumption in Corollary \ref{aban-prob} and the scaling assumption $(\mathbb A3)$.
\end{enumerate}
\end{example}
\begin{remark}
Corollary \ref{aban-prob} generalizes Propositions 2 and 3 in \cite{LeeWard2019POMS}, which are established under the Poisson arrivals (i.e., $M/GI/1+GI$ model) and essentially follow from the results in \cite{HuangGurvich18}.  
\end{remark}

\section{{Uniform Moment Estimates}} \label{section:UniformMoment}

We use the subscript $x$ to denote the scaled Markov process $\wt{\mathbb X}^n$ has an initial state  $(\wt\tau^n(0), \wt V^n(0))= (\tau, v)\equiv x\in\mathbb S$. 




\begin{proposition}\label{moment1}
Assume $(\mathbb A1)$--$(\mathbb A2)$.  {Let $q\in[1,p).$}
{There exist a time $t_0\pos$ and a sufficiently large index $n_0$} 
such that for all $t\geq t_0$,
\beq\label{2.1} \lim_{|x|\ra\infty}\sup_{n\ge n_0}\frac{1}{|x|^{{q}}}\E\left[|\wt{\mathbb X}^n_x(t|x|)|^{{q}} \right]=0. 
\eeq
\end{proposition}
\noindent Proposition \ref{moment1} readily yields uniform (in $n$) moment bounds for the scaled process $\wt{\mathbb X}^n$, from which the Lyapunov function methods of Meyn and Tweedie \cite{Meyn:Tweedie:1993[3]} can yield moment bounds uniformly in time $t$. Ultimately, it implies the moment bounds for the stationary distributions uniform in the scaling parameter $n$ and hence the tightness of the collection of the stationary distributions; see Section 5 in Lee {\it et al.}~\cite{LeeWardYe-QueSyst} for details.   


\begin{proof}[\bf {Proofs of Theorems~\ref{stability} and~\ref{convstat}.}] Given the uniform $q$-th moment estimate in Proposition \ref{moment1}, Theorems~\ref{stability} and~\ref{convstat} follow from Section 5 in Lee {\it et al.}~\cite{LeeWardYe-QueSyst} without any modification. 
\end{proof}



{The crux in proving Proposition~\ref{moment1} lies in two versions of pathwise stability results
(Lemmas \ref{lem-10a} and \ref{lem-10b} below), whose intuitive ideas are provided right after stating those results.} 

\noindent\textbf{Proof of Proposition \ref{moment1}}.



\noindent It is sufficient to show that there exists  $t_0 \in \R_+$ such that for all $t \geq t_0$,
\begin{equation} \label{eq:to-show-1}
\lim_{|x|\ra\infty}\sup_{n}\frac{1}{|x|^{{q}}}\E\left[\wt\tau^n_x(t|x|)^{{q}} \right]=0
\end{equation}
and
\begin{equation} \label{eq:to-show-2}
\lim_{|x|\ra\infty}\sup_{n} \frac{1}{|x|^{{q}}}\E\left[\wt{V}^n_x(t|x|)^{{q}} \right]=0.
\end{equation}


{First, we note that the proof of \eqref{eq:to-show-1} follows from the proof of (11) in Lee {\it et al.}~\cite{LeeWardYe-QueSyst} without any modification.}
To complete the proof, we must show (\ref{eq:to-show-2}),
which is more involved than (\ref{eq:to-show-1}),
and proceeds following the approach of
Ye and Yao (\cite{YeYao12-interchange}, Lemma 10 and Proposition 11).
First, we establish two versions of pathwise stability results
(Lemmas \ref{lem-10a} and \ref{lem-10b}),
one for any (fixed) $n$-th system and the other for the whole sequence.
With the moment condition ($\mathbb A$1) on the system primitives,
the pathwise stability results are then turned into the moment stability
in Lemma \ref{lem-10cd}, which finally leads to (\ref{eq:to-show-2}).

\begin{lemma} \label{lem-10a}
\textsc{(Stability of $\wt V^n(\cdot)$ for large (fixed) $n$)} Let $\{r_i \}_{i\geq1}$ be a sequence of numbers such that
$r_i \rightarrow \infty$ as $i \rightarrow \infty$
and assume the sequence of initial states $\{x^i \in {\mathbb{S}}\}_{i\geq1}$
satisfies $|x^i| \le r_i $ for all $i$.
Then, {there exists some $\bar t_0 >0$ such that for any sufficiently large $n$,}
the following holds (with probability one),
\begin{eqnarray}
 \lim_{i\rightarrow \infty}  \frac{1}{r_i} \wt V_{x^i}^n(r_i t) = 0 ,
  \mbox{~u.o.c.}  \mbox{~~for~~} t\ge \bar t_0 . 
   \label{eq-10a-conv}
\end{eqnarray}
\end{lemma}

\begin{lemma} \label{lem-10b}
\textsc{(Stability of $\wt V^n(\cdot)$)}  Let $\{r_n \}$ be a sequence of numbers such that
 $ r_n  \rightarrow \infty$ as $n\rightarrow \infty$
and assume that the sequence of initial states
$\{x^n\in {\mathbb{S}}\}_{n\geq1}$
satisfies $|x^n|\le r_n$.
 Then, {for some $\bar t_0 >0$},
the following holds (with probability one),
\begin{eqnarray}\label{lemma2-conv}
 \lim_{n\to\infty}\frac{1}{r_n} \wt V_{x^n}^n(r_n t)= 0,
  \mbox{~u.o.c. for ~}  t\ge \bar t_0.
\end{eqnarray}
\end{lemma}

\begin{lemma} \label{lem-10cd}
\textsc{(Moment stability)}  
{The following conclusions hold for some $\bar t_0>0$.}
\\
(a) Letting $\{r_i\}$ and $\{x^i\}$ as in Lemma \ref{lem-10a},
{and the index $n$ be sufficiently large,}
\begin{eqnarray}
   \lim_{i\rightarrow \infty}
     \E  \frac{1}{r_i^q}
   \widetilde V_{x^i}^n(r_i t)^q =0 ,
    \mbox{~~for~~}  t\ge \bar t_0 . 
    \label{eq-m-stab-i}
\end{eqnarray}
(b) Letting $\{r_n \}$ and $\{x^n\}$ as in Lemma \ref{lem-10b},
\begin{eqnarray}
   \lim_{n\rightarrow \infty}
     \E \frac{1}{r_n^q}
   \widetilde V_{x^n}^n(r_n t)^q =0 ,
    \mbox{~~for~~} t\geq \bar t_0 . 
     \label{eq-m-stab-n}
\end{eqnarray}
\end{lemma}

While the proofs of the above three lemmas are provided
in Section \ref{section:LemmaProof},
we provide some intuitions here, for Lemmas \ref{lem-10a} and \ref{lem-10b}
in particular.
Consider the $n$-th (original) system, $V^n(t)$.
The exogeneous arrival rate of customers is increased by $n$ times to $ n \lambda$
(cf. the equation (\ref{254})).
When the workload is $V^n(t) = v$, the arriving customer will abandon with 
a probability $F^n(v)$, and the ``effective'' arrival rate {\it of customers}
becomes $n\lambda (1-F^n(v))$. 
Hence, assuming the server is busy with the service rate $\mu^n$, 
the netput rate of the system in terms of \emph{customers}, 
will be $ n\lambda (1-F^n(v)) - \mu^n $.
As the service rate is $\mu^n$ ($\approx n\lambda$), the 
workload (i.e., the required service time) embodied in each customer
is approximately $(\mu^n)^{-1}$, 
and therefore the netput rate of the system in terms of {\it workload}
will become
  $ (\mu^n)^{-1} [n\lambda (1-F^n(v)) - \mu^n] $.
Since the diffusion-scaled workload inflates the original one by $\sqrt{n}$
(i.e., $\wt V^n(t) = \sqrt{n} V^n(t)$), 
the above rate immediately translates to the netput rate 
in terms of diffusion-scaled workload when $\wt V^n(t) = \wt v (=\sqrt{n}v)$:
\beq\label{com} \sqrt{n} (\mu^n)^{-1} [n\lambda (1-F^n(v)) - \mu^n]
  = n(\mu^n)^{-1} [\sqrt{n} (\lambda-\frac{\mu^n}{n})
         - \sqrt{n}F^n(\frac{\wt v}{\sqrt{n}})] .
\eeq
Observe that $\sqrt{n} (\lambda-{\mu^n}/{n})$ approaches $\theta$ 
from ($\mathbb A$2), 
and that $\sqrt{n} F^n({\wt v}/{\sqrt{n}})$ approaches $\lambda H(\wt v)$
from $(\mathbb A3)$ 
and will be (strictly) greater than $\theta$ for sufficiently large $n$ and 
large workload $\tilde v$ from ($\mathbb A4$).
Then, the above rate will be negative.
Consequently, for large $n$, if the diffusion-scaled workload 
starts from a large (scaled) workload state $\wt V^n(0) = r_i$
as specified in Lemma \ref{lem-10a}
(ignoring the residuals for simplicity), it will have a negative drift,
which is indeed below a negative constant:
  $$ n(\mu^n)^{-1} [\sqrt{n} (\lambda-\frac{\mu^n}{n})
         - \sqrt{n}F^n(\frac{r_i}{\sqrt{n}})]
        \le -\textcolor{black}{\bar\sigma} <0.
    $$
This workload, $\wt V^n(t)$, will reach
the ``normal'' operating state after the initial period
with an order of $r_i / \textcolor{black}{\bar\sigma}$.
The normal operating state,
scaled by $1/r_i$ (where $r_i \rightarrow \infty$), will be approximately zero,
and this is characterized by the convergence in \eqref{eq-10a-conv}
in Lemma \ref{lem-10a}.
Simimlar observation applies to Lemma \ref{lem-10b} as well. 
Lastly, Lemma \ref{lem-10cd} plays a pivotal role in establishing the key moment estimate in \eqref{eq:to-show-2}. 
Given  Lemma \ref{lem-10cd}, the proof of \eqref{eq:to-show-2}
repeats the one for Proposition 11 of \cite{YeYao12-interchange}. 
\begin{remark}\label{rem22}
{From a technical standpoint, the assumptions $(\mathbb A2)$ and $(\mathbb A4)$
can be relaxed such that the above drift \eqref{com} 
 converges to a negative constant as both $n$ and $\tilde v$ go to $+\infty$.}  
\end{remark}


%
%

\begin{proof}[\bf Proof of \eqref{eq:to-show-2}]
{Let $\bar t_0$ be the time given in Lemma \ref{lem-10cd}, 
pick any time $t> \bar t_0$,
{and restrict to all sufficiently large $n$.} 
The rest of the proof is identical to the proof of (12) in 
Lee {\it et al.}~\cite{LeeWardYe-QueSyst}.}
\end{proof}

\section{Proofs of Lemmas and Corollary} \label{section:LemmaProof}
{Lemmas~\ref{lem-10a}, \ref{lem-10b} and \ref{lem-10cd} are almost the same as those in Lee {\it et al.}~\cite{LeeWardYe-QueSyst}. The differences are that Lemma~\ref{lem-10a} and Lemma~\ref{lem-10cd}(a) require the index $n$ be sufficiently large and also the presence of time $\bar t_0>0$ (rather than $c/\sqrt n$ and $ \epsilon$ in Lee {\it et al.}~\cite{LeeWardYe-QueSyst}).  The proof of Lemma~\ref{lem-10cd} is identical, and the proof of Lemma~\ref{lem-10a} follows a similar outline but needs modifications to accommodate the patience time scaling assumption $(\mathbb A3)$ and the stability condition $(\mathbb A4)$.  The proof of Lemma~\ref{lem-10b} is where the main difference lies.}

\begin{proof}[\bf {Proof of Lemma~\ref{lem-10a}.}]
{The proof is a modification of the one for Lemma 1 of 
Lee {\it et al.}~\cite{LeeWardYe-QueSyst}.}
We first show that the workload process $\widetilde V^n$, scaled by $r_i$, converges to a fluid limit process. Then, we analyze the fluid limit to reach our conclusion in (\ref{eq-10a-conv}). The first part of this proof basically repeats the corresponding part of Lemma 1 of 
Lee {\it et al.}~\cite{LeeWardYe-QueSyst}. For the sake of completeness, we include the full proof.

{\textbf{Part 1.}} Without loss of generality, assume that as $i \rightarrow \infty$,
$x^i / r_i \rightarrow {\bar x}\equiv (\bar \tau, \bar v(0))$
with $|{\bar x}| = \bar \tau + \bar v(0) \le 1$; otherwise, it suffices
to consider any  convergent subsequence. 
Fix the index $n$ throughout the proof.
We also omit the index $n$ and the subscript $x^i$
whenever it does not cause any confusion.
For the $i$-th copy of the $n$-th system,
write the offered waiting time as:
\begin{eqnarray}
 \frac{1}{r_i} \wt V_{x^i}^n(r_i t)  & \equiv &
    \wt v_i(t) = \phi_i(t) + \eta_i(t),
    \mbox{~~~~with}
   \label{eq-lee-ward-8} \\
 \phi_i(t) &  \equiv &
   \frac{v^i}{r_i}
   + \frac{n\sqrt{n}}{\mu^n} \cdot \frac{1}{r_i n}
     \sum_{j=1}^{A^n(r_i t)} (1- F^n(V^n(t_j^n -)) )
   - \sqrt{n} t
     \nonumber \\
 && + \frac{n\sqrt{n}}{\mu^n} \cdot \frac{1}{r_i} \left(
       S^n(A^n(r_i t)) - S_{d}^n(A^n(r_i t))
       - \frac{1}{n} M^n(A^n(r_i t))
     \right),
   \label{eq-lee-ward-8b} \\
 \eta_i(t) &  \equiv &
   \frac{\sqrt{n}}{r_i} \int_0^{r_i t} {\bf 1}_{\{V^n(s)=0\}} ds .
   \label{eq-lee-ward-8c}
\end{eqnarray}
In \eqref{eq-lee-ward-8b}, the centered processes $\{S^n(\cdot)\}$, $\{S^n_d(\cdot)\}$, $\{M^n(\cdot)\}$ are defined as follows. For $k\in \N$,
 \[ S^n(k)\equiv {\frac1n\sig{j=1}{k}(v_j-1)},\quad S^n_{d}(k)\equiv \frac1n {\sig{j=1}{k}(v_j-1)}\mathbf{1}_{\{V^n(t^n_j-)\geq d_j^n\}},
   \]
\[
 M^n(k)=\sig{j=1}{k}\left[\mathbf{1}_{[V^n(t^n_j-)\geq d_j^n]}-F^n(V^n(t^n_j-))\right]. 
\]

First, estimate the item associated with the arrival in the above
(in the first summation):
\begin{eqnarray}
 \frac{A^n(r_i t)}{r_i n}
   &=& \frac{1}{r_i n} \left( A^n(r_i t)
     - \lambda n (r_i t - \frac{\tau_i}{\sqrt{n}} \wedge r_i t )
     \right)
     + \lambda (t - \frac{\tau_i}{r_i \sqrt{n}} \wedge t )
     \nonumber \\
 && \rightarrow
   \lambda (t - \frac{\bar \tau}{\sqrt{n}} \wedge t ),
   \mbox{~~~~as } i \rightarrow \infty
   \mbox{~~a.s.}
     \label{eq-arr-fluid}
\end{eqnarray}

Second, denote the term associated with the arrival and abandonment as
\begin{eqnarray}
  && \xi_i(t)
    \equiv \frac{1}{r_i n} \sum_{j=1}^{A^n(r_i t)} (1- F^n(V^n(t_j^n -))).
   \label{eq-xibar}
\end{eqnarray}
Observe that for any $0\le t_1 < t_2 $, we have
\begin{eqnarray}
  &&  0 \le \xi_i(t_2) - \xi_i(t_1)
   \le \frac{1}{r_i n} ( A^n(r_i t_2) - A^n(r_i t_1) ).
    \label{eq-xi-lowbdd}
\end{eqnarray}
From (\ref{eq-arr-fluid}), we note that the right-hand side in
the above converges uniformly to
 $\lambda (t_2 - t_1
  - \frac{\bar \tau}{\sqrt{n}} \wedge t_2
  + \frac{\bar \tau}{\sqrt{n}} \wedge t_1 )$.
Therefore, any subsequence of $i$ contains a further subsequence such
that as $i\rightarrow \infty$ along the further subsequence,
we have the weak convergence
\begin{eqnarray}
  && \xi_i(\cdot) \Rightarrow {\bar \xi}(\cdot) , \,\,\mbox{ in }\,\, D(\R) \,\,\mbox{ as } \,\, i\rightarrow \infty,
    \nonumber
\end{eqnarray}
where the limit ${\bar \xi}(\cdot)$ is Lipschitz continuous (recall \eqref{eq-xi-lowbdd})
 with a Lipschitz constant $\lambda$ (with probability one).
Without loss of generality, we can assume the above convergence is along the full sequence,
and furthermore, by using the coupling technique, we can further assume
the convergence is almost surely:
\begin{eqnarray}
  && \xi_i(t) \rightarrow {\bar \xi}(t) ,
    \mbox{~~~~as } i\rightarrow \infty
    \mbox{~~~a.s}.
    \nonumber
\end{eqnarray}

Third, for the martingale terms, we have
as $i\rightarrow \infty$ with probability one,
\begin{eqnarray}
 && \frac{1}{r_i } \left(
       S^n(A^n(r_i t)) - S_{d}^n(A^n(r_i t))
       - \frac{1}{n} M^n(A^n(r_i t)
     \right)
     \rightarrow 0 .
\end{eqnarray}

Putting the above convergences together yields,
as $i\rightarrow \infty$,
\begin{eqnarray}
 && \phi_i(t) \rightarrow  \bar\phi(t)  \equiv
   \bar v(0)
   + \frac{n\sqrt{n}}{\mu^n} {\bar \xi}(t)
   - \sqrt{n} t ,
   \mbox{~~~u.o.c. of } \,\,t\ge0.
     \label{eq-xi-bar-conv}
\end{eqnarray}
Note from \eqref{eq-lee-ward-8}--\eqref{eq-lee-ward-8c} that
the tuple $( \wt v_i(t), \phi_i(t), \eta_i(t))_{t\geq0}$ satisfies the
one-dimensional linear Skorokhod problem (cf. \S6.2 of \cite{ChenYao01}):
\begin{eqnarray}
 && \wt v_i(t) = \phi_i(t) + \eta_i(t) \ge0, ~~
    d \eta_i(t) \ge0
      \mbox{~with~ } \eta_i(0)=0,  ~~
  \wt v_i(t) d \eta_i(t) = 0 .
    \nonumber
\end{eqnarray}
Hence, by invoking the Lipschitz continuity of the Skorokhod mapping
(cf. Theorem 6.1 of \cite{ChenYao01}),
the convergence in (\ref{eq-xi-bar-conv}) implies
\begin{eqnarray}
 && \frac{1}{r_i} \wt V_{x^i}^n(r_i t)
   \rightarrow \bar v(t)
   \mbox{~~and~~}
   \eta_i(t) \rightarrow  \bar\eta(t)
   \mbox{~~~u.o.c. of }\,\, t\ge0,
     \label{eq-V-fluidlimit}
\end{eqnarray}
with the limit satisfying the Skorokhod problem as well:
\begin{eqnarray}
 && \bar v(t) = \bar\phi(t) + \bar\eta(t) \ge0, ~~
   d \bar\eta(t) \ge0
     \mbox{~with~} \bar\eta(0) =0, ~~
   \bar v(t) d \bar\eta(t) = 0 .
     \label{eq-sko_fluid}
\end{eqnarray}

Next, we further examine the limit $\bar \xi(\cdot)$
following the approach of Chen and Ye
(\cite{ChenYe12}, Proposition 3(b)).
From (\ref{eq-arr-fluid}) and (\ref{eq-xibar}),
and noting that
 $ \xi_i(t) \le \frac{A^n(r_i t)}{r_i n}$,
we have
\begin{eqnarray}
  \bar \xi(t) = 0, ~~ 0\le t \le \frac{\bar \tau}{\sqrt{n}} .
   \label{eq-xi-inittime}
\end{eqnarray}
Now, consider any regular time $t_1 > {\bar \tau}/{\sqrt{n}}$,
at which all processes concerned, i.e.,
$\bar v(\cdot), \bar\phi(\cdot), $ and $\bar\eta(\cdot)$ are differentiable,
and $\bar v(t_1) >0$.
Note that the Lipschitz continuity of $\bar \xi_i(\cdot)$
implies that $(\bar v(\cdot), \bar\phi(\cdot), \bar\eta(\cdot))$ are also
Lipschitz continuous.
Therefore, we can find (small) constants $\epsilon>0$ and $\delta>0$
such that the following inequality holds for all sufficiently large
$i$:
\begin{eqnarray}
 \wt v_i(t_2) > \epsilon
   \mbox{~~i.e., }
   \wt V^n(r_i t_2) > r_i \epsilon ,
   ~~~~ t_2 \in [t_1, t_1 + \delta) .
    \label{eq-vi-epsilon}
\end{eqnarray}
Observe that if the $j$-th arrival falls between
$A^n(r_i t_1)+1$ and  $A^n(r_i t_2)$,
then its arrival time, $t_j^n$, shall also falls between the
corresponding time epochs, i.e.,
$r_i t_1 < t_j^n \le r_i t_2$.
Given the estimate in (\ref{eq-vi-epsilon}),
this implies the following estimate holds:
  $$\wt V^n(t_j^n) > r_i \epsilon .$$
Consequently, we have for all sufficiently large $i$ that
\begin{eqnarray}
   \xi_i(t_2) - \xi_i(t_1)
    &=& \frac{1}{r_i n} \sum_{j= A^n(r_i t_1)+1}^{A^n(r_i t_2)}
    (1 - F^n(V^n(t_j^n -)) )
     \nonumber \\
    &\le&  \frac{1}{r_i n} (A^n(r_i t_2)- A^n(r_i t_1) )
   ( 1 - F^n(r_i \epsilon) ),
   ~~~ t_2 \in [t_1, t_1 + \delta) .
    \nonumber
\end{eqnarray}

{\textbf{Part 2.}}
{From the assumptions in $(\mathbb A3)$ and ($\mathbb A4$), we can find
a constant $\textcolor{black}{\bar\sigma} >0$ such that for any sufficiently large index $n$ and
$x$,
  $$ \sqrt{n} F^n(\frac{x}{\sqrt{n}}) 
    \ge \frac{1}{\lambda}(\theta + 2 \textcolor{black}{\bar\sigma}).
    $$
Putting the above two together yields, for any time $t_2 \in [t_1, t_1 + \delta)$,
\begin{eqnarray}
  \xi_i(t_2) - \xi_i(t_1)
    &\le&  \frac{1}{r_i n} (A^n(r_i t_2)- A^n(r_i t_1) )
   ( 1 - \frac{1}{\sqrt{n} \lambda}(\theta + 2 \textcolor{black}{\bar\sigma}) ) .
    \nonumber
\end{eqnarray}
Taking $i\rightarrow \infty$, this gives
\begin{eqnarray}
 \bar \xi(t_2) - \bar \xi(t_1)
    \le \lambda (t_2 - t_1) ( 1 - \frac{1}{\sqrt{n} \lambda}(\theta + 2 \textcolor{black}{\bar\sigma}) ) .
    \nonumber
\end{eqnarray}
In summary, the above implies
for any regular time $t > {\bar \tau}/{\sqrt{n}}$
with $\bar v(t) >0$,
\begin{eqnarray}
 \frac{d {\bar \xi}(t)}{dt} 
   \le \lambda ( 1 - \frac{1}{\sqrt{n} \lambda}(\theta + 2 \textcolor{black}{\bar\sigma}) ) .
   \label{eq-xi-deri}
\end{eqnarray}
}

{Now, from the properties in
(\ref{eq-xi-bar-conv}), (\ref{eq-sko_fluid}) and (\ref{eq-xi-deri}),
we can see that if $\bar v(t) >0 $ at any time $t\ge \bar\tau /\sqrt{n}$,
\begin{eqnarray}
 \frac{d {\bar v}(t)}{dt} 
   &\le& \frac{n\sqrt{n}}{\mu^n} \lambda 
        ( 1 - \frac{1}{\sqrt{n} \lambda}(\theta + 2 \textcolor{black}{\bar\sigma}) ) 
        - \sqrt{n} 
        \nonumber \\
 &=& \frac{n}{\mu^n} 
       \left( \sqrt{n}(\lambda -\frac{\mu^n}{n}) -\theta - 2\lambda\textcolor{black}{\bar\sigma}) \right) 
   \nonumber 
\end{eqnarray}
According to the assumption $(\mathbb A3)$, if $n$ is sufficienly large, we
can ensure $\sqrt{n}(\lambda -\frac{\mu^n}{n}) -\theta \le \lambda\textcolor{black}{\bar\sigma} $.
Therefore, the above implies
\begin{eqnarray}
 \frac{d {\bar v}(t)}{dt} \le - \textcolor{black}{\bar\sigma} .
    \label{eq-vt-deri+}
\end{eqnarray}
Moreover, combined with the property (\ref{eq-xi-inittime}), the above 
actually holds for any time $t\ge0$ (i.e., including the initial transition period)
such that $\bar v(t) >0 $.
}

{Finally, using the properties in (\ref{eq-sko_fluid}) and (\ref{eq-vt-deri+}),
it is direct to show that
 $\bar v(t) = 0 $
for $t \ge \bar v(0)/ \textcolor{black}{\bar\sigma}$.
This property, along with the convergence in (\ref{eq-V-fluidlimit})
yields the desired convergence in (\ref{eq-10a-conv})
with a constant time $\bar t_0 \ge 1/\textcolor{black}{\bar\sigma}$. 
}
\end{proof}

\begin{proof}[\bf {Proof of Lemma~\ref{lem-10b}.}]
Without loss of generality, we assume $\theta \ge 0$ in ($\mathbb A4$); 
otherwise, the $GI/GI/1+GI$  system is dominated by a stable $GI/GI/1$
system that has the model primitives but no customer abandonment.
Our proof is composed of two parts.

\textbf{Part 1.}
Consider a special case satisfying the extra condition:
there exist a sufficiently large $K^* >0$ and 
a possibly small $\textcolor{black}{\bar\sigma} >0$ (with an abuse of notation, in that it was used in the proof of Lemma~\ref{lem-10a}) such that 
the patience time distributions satisfy the followings
for all $x \ge K^*$ and all sufficiently large $n$,
\begin{eqnarray}
 && \sqrt{n} F^n \left( \frac{x}{\sqrt{n}} \right)
    =  \frac{\theta}{\lambda} + \textcolor{black}{\bar\sigma}.
    \label{assume_specialcase}
\end{eqnarray}
For ease of presentation, we assume the above holds for all $n$.
Below, we prove the lemma under this extra assumption in two steps.

{\it Step 1.}
For the $n$-th system, 
write the offered waiting time as:
\begin{eqnarray}
 \frac{1}{r_n} \wt V_{x^n}^n(r_n t)  & \equiv &
    \wt v_n(t) = \phi_n(t) + \eta_n(t),
    \mbox{~~~~with}
   \label{eq-lee-ward-88} \\
 \phi_n(t) &  = &   \frac{v^n}{r_n} 
   + \phi_{n,1}(t) - \phi_{n,2}(t) + \phi_{n,3}(t) ,
     \nonumber \\
 \phi_{n,1}(t) & = &  
    \frac{n\sqrt{n}}{\mu^n} \cdot \frac{A^n(r_n t)}{r_n n}
    - \sqrt{n} t ,
     \nonumber \\
 \phi_{n,2}(t) & = &  
     \frac{n\sqrt{n}}{\mu^n} \cdot \frac{1}{r_n n}
     \sum_{j=1}^{A^n(r_n t)}  F^n(V^n(t_j^n -)) ,
     \label{phi2-express} \\
 \phi_{n,3}(t) & = & 
    \frac{n\sqrt{n}}{\mu^n} \cdot \frac{1}{r_n} \left(
       S^n(A^n(r_n t)) - S_{d}^n(A^n(r_n t))
       - \frac{1}{n} M^n(A^n(r_n t))
     \right),
   \label{eq-lee-ward-8bb} \\
 \eta_n(t) &  = &
   \frac{\sqrt{n}}{r_n} \int_0^{r_n t} {\bf 1}_{\{V^n(s)=0\}} ds .
   \label{eq-lee-ward-8cc}
\end{eqnarray}

Next, we inspect the terms in $\phi_n(t)$. 
For the initial (diffusion-scaled) states, we assume without loss of generality that
\begin{eqnarray}
 && \frac{x^n}{r_n}= \frac{(\nu^n,\tau^n)}{r_n}
    \rightarrow (\bar \nu, \bar \tau). 
      \label{initialstate_converge}
\end{eqnarray}
Otherwise, we can consider a convergent subsequent.
We can also see that the last term converges to 0, i.e.,
\begin{eqnarray}
 && \phi_{n,3}(t) \rightarrow 0,
      \mbox{ u.o.c. of } t\ge 0. 
    \label{phi_n_3_converg}
\end{eqnarray}

For the term $\phi_{n,1}(t)$, similar to the display (37) of 
Lee {\it et al.}~\cite{LeeWardYe-QueSyst},
we have
\begin{eqnarray}
 \frac{A^n(r_n t)}{r_n n}
   &=& \frac{1}{r_n n} \left( A^n(r_n t)
     - \lambda n (r_n t - \frac{\tau^n}{\sqrt{n}} \wedge r_n t )
     \right)
     + \lambda (t - \frac{\tau^n}{r_n \sqrt{n}} \wedge t ).
     \nonumber 
\end{eqnarray}
Using the above, we can write $\phi_{n,1}(t)$ as,
\begin{eqnarray}
 \phi_{n,1}(t)
   &=& \frac{n}{\mu^n} \frac{1}{\sqrt{r_n}} \cdot
     \left[ \frac{1}{\sqrt{r_n n}} \left( A^n(r_n t)
     - \lambda n (r_n t - \frac{\tau^n}{\sqrt{n}} \wedge r_n t )
     \right) \right]
     \nonumber \\
   &&  + \frac{n\lambda}{\mu^n} \cdot \sqrt{n}
          (1- \frac{\mu^n}{n\lambda})t
     - \frac{n\lambda}{\mu^n} ( \frac{\tau^n}{r_n } \wedge \sqrt{n}t ).
     \nonumber 
\end{eqnarray}
According to the functional central limit theorem, the term in the 
squared parentheses converge to a Brownian motion weakly. 
Therefore, as $r_n\rightarrow \infty$, we have
\begin{eqnarray}
   &&  \frac{1}{\sqrt{r_n}} \cdot
     \left[ \frac{1}{\sqrt{r_n n}} \left( A^n(r_n t)
     - \lambda n (r_n t - \frac{\tau^n}{\sqrt{n}} \wedge r_n t )
     \right) \right]
      \rightarrow 0, 
      \mbox{ u.o.c. of } t\ge 0. 
     \nonumber 
\end{eqnarray}
From the heavy traffic condition, we have
\begin{eqnarray}
   &&  \sqrt{n} (1- \frac{\mu^n}{n\lambda})
     \rightarrow 
     \frac{\theta}{\lambda} .
     \nonumber 
\end{eqnarray}
For the last term, we have for $t=0$,
\begin{eqnarray}
   && \frac{\tau^n}{r_n } \wedge \sqrt{n}t = 0,
     \nonumber 
\end{eqnarray}
and for $t>0$,
\begin{eqnarray}
   && \frac{\tau^n}{r_n } \wedge \sqrt{n}t 
      \rightarrow \bar \tau .
     \nonumber 
\end{eqnarray}
In summary, we have for $t=0$,
\begin{eqnarray}
   && \phi_{n,1}(0) = 0,
     \label{phi_1_0_conv} 
\end{eqnarray}
and for $t>0$,
\begin{eqnarray}
   && \phi_{n,1}(t) 
     \rightarrow \frac{\theta}{\lambda} t - \bar \tau ,
           \mbox{ u.o.c.} 
     \label{phi_1_t_conv} 
\end{eqnarray}

For the term $\phi_{n,2}(t)$, note that for $0\le t_1 < t_2$,
\begin{eqnarray}
 \phi_{n,2}(t_2) - \phi_{n,2}(t_1) & = & 
     \frac{n}{\mu^n} \cdot \frac{1}{r_n n}
     \sum_{j= A^n(r_n t_1)+1}^{A^n(r_n t_2)} \sqrt{n} F^n(V^n(t_j^n -)) 
     \nonumber \\
 & \le & \frac{n}{\mu^n} \cdot \frac{1}{r_n n}
     \left( A^n(r_n t_2) - A^n(r_n t_1) \right) 
      \left( \frac{\theta}{\lambda} + \textcolor{black}{\bar\sigma} \right)
       \nonumber \\
 & \rightarrow & 
    (t_2 - t_1) \left( \frac{\theta}{\lambda} + \textcolor{black}{\bar\sigma} \right)
\end{eqnarray}
As $\phi_{n,2}(t)$ is increasing in $t$, the above implies
\begin{eqnarray}
 && \phi_{n,2}(t) \rightarrow \hat \phi_2(t),
      \mbox{ u.o.c. of } t\ge 0. 
    \label{converg_to_hatphi2}
\end{eqnarray}
where the limit $\hat \phi_2(t)$ is increasing 
and Lipschitz continuous in $t\ge0$, with $\hat \phi_{2}(0)=0$.

Putting (\ref{initialstate_converge}, \ref{phi_n_3_converg},
\ref{phi_1_0_conv}, \ref{phi_1_t_conv}, \ref{converg_to_hatphi2})
together, we have
\begin{eqnarray}
 && \phi_{n}(t) \rightarrow \hat \phi(t),
      \mbox{ u.o.c. of } t > 0, 
    \label{conv_to_hatphi_t}
\end{eqnarray}
where the limit $\hat \phi(t)$, $t\ge0$, is Lipschitz continuous in
$(0, \infty)$, and satisfies
\begin{eqnarray}
 && \hat\phi(0) = \bar \nu, ~~~
    \hat\phi(0+) = \max(0, \bar\nu -\bar\tau). 
    \label{conv_to_hatphi_0}
\end{eqnarray}
The above combined with the reflecting mapping given in (4),
we have 
\begin{eqnarray}
 && (\tilde v_n(t), \phi_{n}(t), \eta_{n}(t)) 
    \rightarrow 
     (\hat v(t), \hat \phi(t), \hat\eta(t)),
      \mbox{ u.o.c. of } t > 0, 
    \label{limit-sko}
\end{eqnarray}
where the limit forms a standard one-dimensional Skorohod mapping.

{\it Step 2.}
To prove the Lemma for the special case, it suffices to
show the followings for any $t>0$:
\begin{eqnarray}
 && {\frac{d}{dt} } \hat \phi(t) = -\textcolor{black}{\bar\sigma}, 
      \mbox{ ~~when } \hat v(t) > 0. 
    \label{negativeslop}
\end{eqnarray}
{Indeed, the above property, along with
the characterization in (\ref{conv_to_hatphi_t}-\ref{limit-sko}), 
imply 
$\hat v(t_0 + \cdot)$ for some $t_0>0$,
which then justifies the conclusion in (\ref{lemma2-conv}).
}

Fix any $t>0$ such that $\hat v(t) > 0$.
Then, we can find some (small) interval $[t_1, t_2]$
such that $\hat v(t') > 2 \epsilon$
for any $t' \in [t_1, t_2]$ and for some constant $\epsilon>0$.
Hence, for sufficiently large $n$, we have
for all $t' \in [t_1, t_2]$
  $$\frac{1}{r_n} \wt V^n(r_n t') > \epsilon .$$
Note, for any $t_j^n \in [r_n t_1, r_n t_2]$
(i.e., $t_j^n /r_n\in [t_1, t_2]$), we have for sufficiently large $n$,
  $$ V^n(t_j^n -) = \frac{1}{\sqrt{n}} \tilde V^n(t_j^n-)
    \ge \frac{r_n \epsilon}{\sqrt{n}} 
     \ge \frac{K^*}{\sqrt{n}} . $$
Therefore, from our assumption (\ref{assume_specialcase}),
we have for sufficiently large $n$,
\begin{eqnarray}
 && \sqrt{n} F^n(V^n(t_j^n-)) 
   = \frac{\theta}{\lambda} + \textcolor{black}{\bar\sigma}  . 
    \nonumber
\end{eqnarray}
From (\ref{phi2-express}, \ref{converg_to_hatphi2}), we have 
\begin{equation}
   \frac{d}{dt} \hat\phi_2(t) = \frac{\theta}{\lambda} + \textcolor{black}{\bar\sigma},
    ~~~ t>0 .
    \nonumber
\end{equation}
The above equality and the convergence (\ref{phi_1_t_conv})
together yield (\ref{negativeslop}) immediately.



\textbf{Part 2.}
To prove the lemma for the general case, 
we construct dominating systems that satisfy the condition in 
the special case first. The dominating systems have the same settings as the original $GI/GI/1$ queues 
under study, except the patience times. 
That is, they are driven by the same sequences of arrival times $\{ t_j^n \}$ 
and service times $\{ v_j^n \}$.
However, the patience times in the dominating systems,
denoted as $d_j^{n*}$ for the $j$-th arrival of the $n$-th queue,
are generated from 
those in the original queues using the inverse transformation method 
as follows. 

First, 
under the stability condition ($\mathbb A4$)
and the model assumption $(\mathbb A3)$, we can choose 
a sufficiently large $K^* >0$ and 
a possibly small $\textcolor{black}{\bar\sigma} >0$ such that 
the original patience time distributions satisfy the followings
for  all sufficiently large $n$,
\begin{eqnarray}
 && \sqrt{n} F^n \left( \frac{x}{\sqrt{n}} \right)
    \ge  \frac{\theta}{\lambda} + \textcolor{black}{\bar\sigma} ,
    \mbox{~~~for } x \ge K^* .
    \label{assume_specialcase_2}
\end{eqnarray}
Now, we specify the distribution $F^{n*}(\cdot)$ 
of patience times for the $n$-th
dominating queue as 
\begin{eqnarray}
 && \sqrt{n} F^{n*} \left( \frac{x}{\sqrt{n}} \right)
     = \min \left\{ \sqrt{n} F^n \left( \frac{x}{\sqrt{n}} \right) ,
         \frac{\theta}{\lambda} + \textcolor{black}{\bar\sigma}
     \right\} ,
     \mbox{~~~ for } x \ge0 .
\end{eqnarray}
Then, the distribution $F^{n*}(\cdot)$ satisfies the condition 
(\ref{assume_specialcase}), i.e.,
\begin{eqnarray}
 && \sqrt{n} F^{n*} \left( \frac{x}{\sqrt{n}} \right)
    =  \frac{\theta}{\lambda} + \textcolor{black}{\bar\sigma}.
    \label{assume_specialcase_3}
\end{eqnarray}
Below, for ease of presentation, we assume the above equality holds
for all $n$.

Next, we construct the patience times in the dominating queues as, 
{for each $\omega\in\Omega, n\geq1, j\geq 1$,}
\begin{eqnarray}
 && d_j^{n*}{(\omega)} = \inf \{x: F^{n*}(x) = F^n(d_j^n{(\omega)})\} .
    \label{patiencetime_domin}
\end{eqnarray}
Note that a patience time is allowed to be infinity, in which case
the customer will never abandon. 
Observe that 
\begin{itemize}
 \item for each $n$, the sequence of patience times
 $\{d_j^{n*}, j=1, 2, \cdots \}$ are i.i.d.~and follow
 the distribution $F^{n*}(\cdot)$, and
 
 \item the dominating systems have longer patience times (for all samples):
 $d_j^{n*} \ge d_j^n$.
\end{itemize}

Now, as the dominating systems have the patience time distribution
$F^{n*}(\cdot)$ that satisfy the property (\ref{assume_specialcase_3}),
they fit into the special case specified in the condition
(\ref{assume_specialcase}). 
Hence, the conclusion in ``modified Lemma 2'' holds for 
the sequence of dominating queues, i.e., 
the associated workloads,
denoted as $V^{n*}(t)$ (or $\tilde V^{n*}(t)$ under the diffusion scaling),
satisfy the convergence property (\ref{lemma2-conv}):
\begin{eqnarray}
 && \lim_{n\to\infty}\frac{1}{r_n} \wt V_{x^n}^{n*}(r_n t)= 0,
   \mbox{~~~u.o.c.~for } t\ge t_0. 
    \label{eq-domin_cov2zero}
\end{eqnarray}
Moreover, 
since each customer in the dominating queues has a longer patience time
than the corresponding one in the original queues,
the workload of each dominating queue must be (roughly) as much as 
the original one.
Indeed, we can show:
\begin{eqnarray}
 && V^{n*}(t) \ge V^n(t) - \nu^{n,\max}(t) ,
   \label{eq-dominating_0613}
\end{eqnarray}
where $\nu^{n,\max}(t)$ is the maximum of all workloads that customers have
brought in by time $t\ge0$,
 $$ \nu^{n,\max}(t) = \max \{ \nu_j^n : j=1, \cdots, A^n(t) \} . $$
From Lemma 5.1 of \cite{Bramson:1998}, we have, with probability 1,
the following u.o.c.~convergence as $n\rightarrow \infty$:
\begin{eqnarray}
 && \frac{\sqrt{n}}{r_n} \nu^{n,\max}(r_n t)
  ~ \left[ = \frac{1}{r_n} \tilde \nu^{n,\max}(r_n t) \right]
   \rightarrow 0 . 
     \label{eq-residual_conv0}
\end{eqnarray}
(Here, $ r_n n $ corresponds to $r^2$ of  Lemma 5.1 Bramson \cite{Bramson:1998}.  In our definition of $\nu^{n,\max}$, 
 $A^n$ refers to the {arrival} (with $\nu_j^n$ being the {service} time), 
 while in Bramson \cite{Bramson:1998} the corresponding variable, $v^{r,T, \,max}$,
 connects to the {service} times and the {service} (renewal) process.) 
Consequently, the conclusion in Lemma \ref{lem-10b} follows from
(\ref{eq-domin_cov2zero}-\ref{eq-residual_conv0}).
\medskip
It remains to prove (\ref{eq-dominating_0613}), and this is done by induction.
First, note that the inequality holds for $t= t_0^n \equiv 0$, as both systems share
the same initial state ($V^n(0) = V^{n*}(0)$).
Next, we show that if for any $J\ge0$, it holds for $t= t_J^n$:
\begin{eqnarray}
 && V^{n*}(t_J^n) \ge V^n(t_J^n) - \nu^{n,\max}(t_J^n) ,
   \label{eq-dominating_0613-t_J}
\end{eqnarray}
then, it can be extended to the time $t\in (t_J^n , t_{J+1}^n]$.

When $t\in (t_J^n , t_{J+1}^n)$, by using the evolution equation \eqref{2.1.1} and the relationship $A^n(t_j^n) = j$
(which hold for both the original and the dominating systems),
we write
\begin{eqnarray}
 V^n(t) & = & V^n(t_J^n) - \int_{t_J^n}^{t} {\bf 1}_{[V^n(s)>0]} ds 
    \nonumber \\
 & = & \left\{
        \begin{tabular}{ll}
           $ V^n(t_J^n) - (t - t_J^n) $, 
             & $ t- t_J^n <  V^n(t_J^n)$ ,\\
           $ 0 $, 
             & $ t- t_J^n \le V^n(t_J^n) $ ,\\
        \end{tabular}
       \right.
   \label{eq-t_J_J+1_a}
\end{eqnarray}
and similarly, 
\begin{eqnarray}
 V^{n*}(t) 
 & = & \left\{
        \begin{tabular}{ll}
           $ V^{n*}(t_J^n) - (t - t_J^n) $, 
             & $ t- t_J^n <  V^{n*}(t_J^n)$ ,\\
           $ 0 $, 
             & $ t- t_J^n \le V^{n*}(t_J^n) $ .\\
        \end{tabular}
       \right.
   \label{eq-t_J_J+1_b}
\end{eqnarray}
By carefully examining the sample paths as characterized in 
(\ref{eq-t_J_J+1_a}) and (\ref{eq-t_J_J+1_b}),
along with the inductive assumption (\ref{eq-dominating_0613-t_J}),
we can see that the inquality (\ref{eq-dominating_0613}) holds 
for $t\in (t_J^n , t_{J+1}^n)$.
Particularly, this gives
\begin{eqnarray}
 && V^{n*}(t_{J+1}^n -) \ge V^n(t_{J+1}^n -) - \nu^{n,\max}(t_{J+1}^n -) .
   \label{eq-dominating_0613-t_J+1_before}
\end{eqnarray}

For $t = t_{J+1}^n$, we apply the evoluation equation again to write,
\begin{eqnarray}
 V^n(t_{J+1}^n)  & = & V^n(t_{J+1}^n -)
   + \nu_{J+1}^n {\bf 1}_{[V^n(t_{J+1}^n -) < d_{J+1}^n]},
   \label{eq-t_J+1_c} \\
 V^{n*}(t_{J+1}^n)  & = & V^{n*}(t_{J+1}^n -)
   + \nu_{J+1}^n {\bf 1}_{[V^{n*}(t_{J+1}^n -) < d_{J+1}^{n*}]}.
   \label{eq-t_J+1_d}
\end{eqnarray}
Given that $d_{J+1}^{n*} \ge d_{J+1}^n$, 
the relationships in 
(\ref{eq-dominating_0613-t_J+1_before}-\ref{eq-t_J+1_d})
implies the inquality (\ref{eq-dominating_0613}) for $t = t_{J+1}^n$.

\end{proof}

\begin{proof}[\bf {Proof of Lemma~\ref{lem-10cd}.}]
The proof is identical to the one for Lemma 3 of 
Lee {\it et al.}~\cite{LeeWardYe-QueSyst}.
\end{proof}

\begin{proof}[\bf {Proof of Corollary~\ref{aban-prob}.}]
{From the convergence in Theorem~\ref{convstat}(a) and a generalized continuous mapping theorem (see, e.g., Theorem 3.4.4 in Whitt~\cite{Whitt2002}), it suffices to (i) establish a uniform integrability of $\{\sqrt n F^n(\wt V^n(\infty)/\sqrt n)\}_{ n\geq1}$ and (ii) verify the convergence $\sqrt n F^n (x_n/\sqrt n) \ra H(x)$ whenever $x_n\ra x$ as $n\ra\infty$ for a nonnegative sequence $\{x_n\}_{n\geq1}$.  The first part (i): The uniform integrability  of $\{\sqrt n F^n(\wt V^n(\infty)/\sqrt n)\}_{ n\geq1}$ follows immediately from the assumed polynomial growth condition $\sqrt n F^n(x/\sqrt n)\leq C (1+x^m)$ and the uniform integrability of  $\{[\wt V^n(\infty)]^m\}_{n\geq1}$ established in Theorem~\ref{convstat}(b). The second part (ii): Since a nonnegative sequence $\{x_n\}_{n\geq1}$ is convergent, there is $M\pos$ such that $x_n\leq M$. Therefore, 
\beqys 
|\sqrt n F^n(x_n/\sqrt n)-H(x)| &\leq& |\sqrt n F^n(x_n/\sqrt n)-H(x_n)| + |H(x_n)-H(x)|\\
&\leq& \sup_{y\in[0,M]}| \sqrt n F^n(y/\sqrt n)-H(y) | +|H(x_n)-H(x)|,  
\eeqys 
where the first term converges to zero due to $(\mathbb A3)$ and the same for the second term from the continuity of $H(\cdot)$ in $(\mathbb A3)$. This completes the proof. 
}
\end{proof}

\noindent{\bf \large Acknowledgements:  }{
Amy Ward is  supported as Charles M. Harper Faculty Fellow at the University of Chicago Booth School of Business,
and Heng-Qing Ye is supported by the HK/RGC Grant 15503519.


\newcommand{\noop}[1]{}

\end{document}